\documentclass[12pt,reqno]{amsart}
\usepackage{amsfonts,amsmath,amscd,amssymb,amsthm}
\usepackage{enumerate,hyperref,perpage,url}
\usepackage[margin=2.2cm, a4paper]{geometry}

\theoremstyle{plain}
\newtheorem{thm}{Theorem}

\newtheorem{cor}[thm]{Corollary}
\newtheorem{lemma}[thm]{Lemma}

\newtheorem{prop}[thm]{Proposition}

\theoremstyle{remark}
\newtheorem*{defn}{\textbf{Definition}}

\numberwithin{equation}{section}

\MakePerPage{footnote} \allowdisplaybreaks 

\newcommand{\C}{\mathbb C}

\newcommand\HH{\mathcal H}

\newcommand\M{\mathbb M}

\newcommand\N{\mathbb N}
\newcommand\NN{\mathcal N}

\newcommand\R{\mathbb R}

\newcommand\T{\mathbb T}

\newcommand{\Vol}{\mathrm{Vol}}

\newcommand{\wt}{\widetilde}

\title[Distribution of the nodal sets of eigenfunctions]{Distribution of the nodal sets of eigenfunctions on analytic manifolds}

\author{Xiaolong Han}
\address{Department of Mathematics, California State University, Northridge, CA 91325, USA}
\email{Xiaolong.Han@csun.edu}

\subjclass[2010]{58J50, 35J05, 35P15}

\keywords{Laplacian eigenfunctions, real zeros, equidistribution of nodal sets}

\thanks{}

\date{}

\begin{document}
\maketitle

\begin{abstract}
The nodal set of the Laplacian eigenfunction has co-dimension one and has finite hypersurface measure on a compact Riemannian manifold. In this paper, we investigate the distribution of the nodal sets of eigenfunctions, when the metric on the manifold is analytic. We prove that if the eigenfunctions are equidistributed at a small scale, then the weak limits of the hypersurface volume form of their nodal sets are comparable to the volume form on the manifold.

In particular, on the negatively curved manifolds with analytic metric and on the tori, we show that in any eigenbasis, there is a full density subsequence of eigenfunctions such that the weak limits of the hypersurface volume form of their nodal sets are comparable to the volume form on the manifold.
\end{abstract}

\section{Introduction}
Let $(\M,g)$ be a compact and smooth Riemannian manifold of dimension $n\ge2$ without boundary. Denote $\Delta=\Delta_g$ the (positive) Laplace-Beltrami operator. Let $\{e_j\}_{j=0}^\infty$ be a sequence of eigenfunctions of $\Delta$ with eigenvalues $0=\lambda_0^2<\lambda_1^2\le\lambda_2^2\le\cdots$, that is, $\Delta e_j=\lambda^2_je_j$. 

Define the nodal set of eigenfunction $e_j$ as
$$\NN_{e_j}:=\{x\in\M:e_j(x)=0\}.$$
Then $\NN_{e_j}$ has co-dimension one in $\M$ and has finite $(n-1)$-$\dim$ hypersurface measure. A well-known problem in the geoemtry of eigenfunctions is to determine the asymptotics of the hypersurface measure of the nodal sets $\NN_{e_j}$ as $\lambda_j\to\infty$. In fact, Yau \cite{Ya} conjectured that
\begin{equation}\label{eq:Yau}
c\lambda_j\le\HH^{n-1}(\NN_{e_j})\le C\lambda_j,
\end{equation}
in which the positive constants $c$ and $C$ depend only on $\M$. Here, $\HH^{n-1}$ is the $(n-1)$-$\dim$ Hausdorff measure.

In this paper, we investigate the distribution of the nodal sets of eigenfunctions on the manifold. In the view of \eqref{eq:Yau}, we need the normalization $1/\lambda_j$ and are interested in the following problem.

Let $dS_{e_j}$ be the $(n-1)$-$\dim$ Riemannian hypersurface volume form on $\NN_{e_j}$. 
\begin{equation}\label{eq:dS}
\text{Determine all the possible weak limits of }\frac1\lambda_jdS_{e_j}\text{ in }\M,
\end{equation}
and their relation with the geometry of $\M$, that is, study the asymptotics 
$$\frac{1}{\lambda_j}\int_{\NN_{e_j}}f\,d\HH^{n-1}\quad\text{for }f\in C(\M)\text{ as }j\to\infty.$$ 
Of course, \eqref{eq:Yau} is reduced to choosing $f\equiv1$ and proving that $\frac{1}{\lambda_j}\int_{\NN_{e_j}}1\,d\HH^{n-1}\approx1$ as $j\to\infty$. So we see that the distribution of the nodal sets in Problem \ref{eq:dS} implies the estimates of their hypersurface measure in Yau's conjecture \eqref{eq:Yau}.

Donnelly-Fefferman \cite{DF} proved \eqref{eq:Yau} if the metric $g$ on $\M$ is analytic. Hence, on an analytic manifold, the sequence in \eqref{eq:dS} has uniformly bounded mass. If the metric $g$ is only smooth on $\M$, Yau's conjecture \eqref{eq:Yau} remains open and only partial results are proven. Therefore, it is not known that the sequence in \eqref{eq:dS} has uniformly bounded mass. We refer to a survey article Zelditch \cite{Z3} for the recent progress and the partial results related to \eqref{eq:Yau} on smooth manifolds.

Motivated by Zelditch \cite[Conjecture 1.7]{Z2}, we study the phenomenon that the nodal sets display equidistribution, more precisely, in Problem \eqref{eq:dS},
\begin{equation}\label{eq:enodal}
\frac{1}{\lambda_j}\int_{\NN_{e_j}}f\,d\HH^{n-1}\to c\int_\M f\,d\Vol\quad\text{for all }f\in C(\M)\text{ as }j\to\infty.
\end{equation}
This would show that the weak limit of the hypersurface volume form of the nodal sets coincides with (possibly a constant multiple of) the volume form on the manifold, which says that the nodal sets tend equidistributed asymptotically on the manifold.

As was pointed in Zelditch \cite{Z2}, equidistribution of the nodal sets \eqref{eq:enodal} is closely related to the global dynamics the geodesic flow on the manifold. It is in general a very difficult problem. However, on an analytic manifold, the eigenfunctions can be extended to a complex neighborhood of the manifold. The problem \eqref{eq:dS} can be simplified if one studies the complex nodal sets of the eigenfunctions in the complex region. In fact, for the quantum ergodic sequence of eigenfunctions $\{u_k\}$, Zelditch \cite[Theorem 1.1 and Corollary 1.2]{Z2} proved an explicit asymptotic of $dS_{u^\C_k}/\lambda_k$, where $u^\C_j$ denotes the complex extension of $u_k$ and $dS_{u^\C_k}$ is the hypersurface volume form of the nodal set of $u^\C_k$ in the complex region. (See the following discussion on quantum ergodic eigenfunctions.)

Zelditch \cite[Conjecture 1.7]{Z2} then proposed that, the real nodal sets of these quantum ergodic eigenfunctions on analytic manifolds should also be equidistributed on $\M$, i.e. \eqref{eq:enodal} is true. Throughout the paper, we assume analyticity of the metric $g$ on $\M$ and address this question involving the equidistribution of the (real) nodal sets of eigenfunctions on $\M$. 

First we need to point out that the precise asymptotic in \eqref{eq:enodal} on $\M$ seems only plausible on some classes of manifolds with specific structure, e.g. arithmetic hyperbolic manifolds, the spheres, the tori, etc.  This is because, if \eqref{eq:enodal} is valid, then choosing $f=1$ implies a precise asymptotic in the hypersurface measure of the nodal sets in \eqref{eq:Yau}, rather than their lower and upper bounds. So in \cite[Conjecture 1.7]{Z2}, the asymptotic in \eqref{eq:enodal} is replaced by a weaker condition as the following comparability.
\begin{equation}\label{eq:compnodal}
\frac{1}{\lambda_j}\int_{\NN_{e_j}}f\,d\HH^{n-1}\approx\int_\M f\,d\Vol\quad\text{for all }f\in C(\M)\text{ as }j\to\infty.
\end{equation}
Our main results in Theorems \ref{thm:nodalinB} and \ref{thm:nodalonM} exactly prove \eqref{eq:compnodal} of eigenfunctions on some manifolds, thus establishing a weaker characterization of the equidistribution phenomenon of their nodal sets on these manifolds.

Now we recall some background on quantum erogidicity that was used in \cite{Z2}: The study of equidistribution of the nodal sets falls in the quantum ergodicity developed by \v Snirel'man \cite{Sn}, Zelditch \cite{Z1}, and Colin de Verdi\`ere \cite{CdV}. That is, on manifolds with ergodic geodesic flow, any eigenbasis $\{e_j\}$ contains a full density subsequence $\{e_{j_k}\}$ which is equidistributed on $\M$, i.e.
\begin{equation}\label{eq:e}
\int_\Omega|e_{j_k}|^2\,d\Vol=\frac{\Vol(\Omega)}{\Vol(\M)}\quad\text{as }k\to\infty
\end{equation}
for all Jordan measurable $\Omega\subset\M$. We call $\{e_{j_k}\}$ a quantum ergodic subsequence of eigenfunctions. Here, We define the density $D$ of a subsequence $J=\{j_k\}\subset\N$ as
$$D(J)=\lim_{N\to\infty}\frac{\#\{j_k<N\}}{N}\quad\text{if it exists}.$$
When $D=1$, we call such subsequence a full density subsequence. Important examples of manifolds with ergodic geodesic flow include negatively curved manifolds, i.e. all the sectional curvatures are negative everywhere.

Quantum ergodicity in fact asserts a stronger result than \eqref{eq:e} that the eigenfunctions $\{e_{j_k}\}$ tend equidistributed in the phase space, roughly speaking, it says that in addition to \eqref{eq:e}, the oscillation of the eigenfunctions also tend equidistributed. (The full statement of which can be defined by microlocal analysis and pseudodifferential operators and we refer to \cite{Sn, Z1, CdV} for details.) Hence, intuitively equidistribution of the eigenfunctions and their oscillation indicate that their nodal set should display some equidistribution phenomenon as in \eqref{eq:compnodal}.

However, because of the singular nature of the nodal sets, i.e. they have co-dimension one on the manifold, \eqref{eq:e} is not sufficient to imply equidistribution of their nodal sets. 

Our main results state that a finer analysis than \eqref{eq:e} implies equidistribution of the nodal sets. Such finer analysis involves equidistribution of eigenfunctions at small scales, that is, \eqref{eq:e} when $\Omega$ are replaced by balls $B(x,r)$ with radius $r\to0$ as the eigenvalues tend to infinity. Because we can control the $L^2$ norm of eigenfunctions in the shrinking balls, we can estimate the nodal sets in the same balls following the argument in Donnelly-Fefferman \cite{DF}, and hence provide estimates of the distribution of the nodal sets at a small scale, which is sufficient for us to prove \eqref{eq:compnodal}. 

Small scale equidistribution has been proved in \cite{LS, Yo, HR1, HR2, H, LR} on some classes of manifolds including arithmetic hyperbolic manifolds,  negatively curved manifolds, the tori, etc. Instead of proving \eqref{eq:compnodal} for the eigenfunctions on these manifolds individually, we present our results in a more general setting, that is, we prove our main theorems assuming small scale equidistribution of a sequence of eigenfunctions. One can then apply the main theorems (Theorems \ref{thm:nodalinB} and \ref{thm:nodalonM}) to the analytic manifolds on which small scale equidistribution holds. 

Let $\{e_{j_k}\}\subset\{e_j\}$ be a subsequence of eigenfunctions. For simplicity, we denote $u_k=e_{j_k}$ and $r_k=r(\lambda_{j_k})$. We first define the small scale functions.
\begin{defn}[Small scale functions]
We say that $r(\lambda):[0,\infty)\to(0,1)$ is a small scale function if $\lambda^{-\rho}\le r(\lambda)<1$ for some $\rho\in(0,1]$ and $r(\lambda)\to0$ as $\lambda\to\infty$.
\end{defn}

Given a small scale function $r(\lambda)$, we define small scale equidistribution of a sequence of eigenfunctions. Denote $B(x,r)$ as a geodesic ball in $\M$ centered at $x$ and with radius $r$.
\begin{defn}[Equidistribution of eigenfunctions at small scales]
Let $\{u_k\}_{k=1}^\infty$ be a sequence of eigenfunctions. We say that $\{u_k\}_{k=1}^\infty$ is equidistributed at a small scale $r(\lambda)$ if
$$\int_{B(x,r_k)}|u_k|^2\,d\Vol=\frac{\Vol(B(x,r_k)}{\Vol(\M)}+o(\Vol(B(x,r_k))\quad\text{as }k\to\infty$$
uniformly for all $x\in\M$.
\end{defn}

However, to prove \eqref{eq:compnodal}, the following uniform comparability condition of $\int_{B(x,r_k)}|u_k|^2$ and $\Vol(B(x,r_k))$ on small balls at a small scale is sufficient for our purpose, that is, for some positive constants $D_1$ and $D_2$ that are independent of $x$ and $k$,
\begin{equation}\label{eq:sse}
D_1\Vol(B(x,r_k))\le\int_{B(x,r_k)}|u_k|^2\,d\Vol\le D_2\Vol(B(x,r_k)\quad\text{as }k\to\infty
\end{equation}
uniformly for all $x\in\M$. 

The main theorem then states that, for a sequence of eigenfunctions that small scale uniform comparability condition in \eqref{eq:sse} is true, then within the balls at a small scale, the hypersurface volume form of their nodal sets are comparable to the volume form on $\M$.
\begin{thm}\label{thm:nodalinB}
On an analytic manifold $\M$, let $\{u_k\}_{k=1}^\infty$ be a sequence of eigenfunctions that is equidistributed at a scale $r(\lambda)$. Then there exist positive constants $E_1$ and $E_2$ depending only on $D_1$ and $D_2$ in \eqref{eq:sse} and $\M$ such that
\begin{equation}\label{eq:nodalinB}
E_1\Vol(B(x,r_k))\le\frac1\lambda_k\HH^{n-1}(\NN_{u_k}\cap B(x,r_k))\le E_2\Vol(B(x,r_k))\quad\text{as }k\to\infty
\end{equation}
uniformly for all $x\in\M$.
\end{thm}

From Theorem \ref{thm:nodalinB}, we derive the second main theorem of this paper.
\begin{thm}\label{thm:nodalonM}
On an analytic manifold $\M$, let $\{u_j\}_{j=1}^\infty$ be a sequence of eigenfunctions that is equidistributed at a small scale $r(\lambda)$. Then there exist positive constants $C_1$ and $C_2$ depending only on $D_1$ and $D_2$ in \eqref{eq:sse} and $\M$ such that
\begin{equation}\label{eq:nodalonM}
C_1\int_\M f\,d\Vol\le\frac1\lambda_k\int_{\NN_{u_k}}f\,d\HH^{n-1}\le C_2\int_\M f\,d\Vol\quad\text{as }k\to\infty
\end{equation}
for any nonnegative function $f\in C(\M)$.
\end{thm}

We now apply Theorems \ref{thm:nodalinB} and \ref{thm:nodalonM} to some classes of manifolds. Han \cite{H} and Hezari-Rivi\`ere \cite{HR1} proved \eqref{eq:compnodal} on all negatively curved manifolds with smooth metric. In particular,
\begin{cor}\label{cor:nodalequivncm}
On a negatively curved manifold $\M$ with analytic metric, let 
$$r(\lambda)=\frac{1}{(\log\lambda)^\alpha}\quad\text{for }\alpha\in\left(0,\frac{1}{2n}\right).$$
There exist positive constants $E_1$, $E_2$, $C_1$, and $C_2$ depending only on $\M$ such that every eigenbasis $\{e_j\}_{j=1}^\infty$ contains a full density subsequence $\{e_{j_k}\}_{k=1}^\infty$ for which 
\begin{itemize}
\item \eqref{eq:nodalinB} is valid for $u_k=e_{j_k}$ and $r_k=r(\lambda_{j_k})$ uniformly for all $x\in\M$;
\item \eqref{eq:nodalonM} is valid for all nonnegative functions $f\in C(\M)$.
\end{itemize}
\end{cor}

On the tori, following Hezari-Rivi\`ere \cite{HR2} and Lester-Rudnick \cite{LR}, we have that
\begin{cor}\label{cor:nodalequivtori}
On the torus $\T^n$ for $n\ge2$, let
$$r(\lambda)=\frac{1}{\lambda^\rho}\quad\text{for }\rho\in\left(0,\frac{1}{n-1}\right).$$
There exist positive constants $E_1$, $E_2$, $C_1$, and $C_2$ depending only on $n$ such that every eigenbasis $\{e_j\}_{j=1}^\infty$ contains a full density subsequence $\{e_{j_k}\}_{k=1}^\infty$ for which 
\begin{itemize}
\item \eqref{eq:nodalinB} is valid for $u_k=e_{j_k}$ and $r_k=r(\lambda_{j_k})$ uniformly for all $x\in\M$;
\item \eqref{eq:nodalonM} is valid for all nonnegative functions $f\in C(\M)$.
\end{itemize}
\end{cor}

Throughout this paper, $A\lesssim B$ ($A\gtrsim B$) means $A\le cB$ ($A\ge cB$) for some constant $c$ depending only on the manifold; $A\approx B$ means $A\lesssim B$ and $B\lesssim A$; the constants $c$ and $C$ may vary from line to line.

\section{Nodal sets of eigenfunctions in balls at small scales}
In this section, we prove Theorem \ref{thm:nodalinB}, that is, if
$$D_1\Vol(B(x,r_k)\le\int_{B(x,r_k)}|u_k|^2\,d\Vol\le D_2\Vol(B(x,r_k),$$
then
$$E_1\Vol(B(x,r_k))\le\frac1\lambda_k\HH^{n-1}(\NN_{u_k}\cap B(x,r_k))\le E_2\Vol(B(x,r_k)),$$
where the positive constants $E_1$ and $E_2$ depend only on $D_1$ and $D_2$ and $\M$. Our argument is to first dilate $u_k$ on $B(x,r_k)$ to a new function on a ball with fixed radius and then to use the estimates of the nodal sets obtained in Donnelly-Fefferman \cite{DF} to the dilated function. Throughout the process, we trace all the dependence on the constants above.

Since $r_k\to0$ as $k\to\infty$, we may assume $B(x,10r_k)$ is contained in one sufficiently small coordinate patch of $\M$, where the metric $g$ can be expanded in power series. For notational convenience, we write $u=u_k$ and $r=r_k$ and identify $x$ as $0$ in this local coordinate patch. Using a dilation, let
$$v(x)=u(rx)\quad\text{for }x\in\{x\in\R^n:|x|\le10\}.$$ 
Then $v$ is an eigenfunction of $\wt\Delta$ in $B(0,10)$ with eigenvalue $\mu^2=r^2\lambda^2$. The metric $\tilde g$ defined on $B(0,10)$ under such a dilation tends to the Euclidean metric as $r\to0$. We can then assume that the absolute values of the section curvatures of $(B(0,10),\tilde g)$ are bounded by $K/2$, where $K$ is an upper bound of the absolute values of the sectional curvatures of $\M$.

Notice that that are positive constants $c_1$ and $c_2$ depending only on $\M$ such that
$$\frac{c_1}{r^n}\int_{B(p,r)}|u|^2\,d\Vol_g\le\int_{B(p,1)}|v|^2\,d\Vol_{\tilde g}\le\frac{c_2}{r^n}\int_{B(p,r)}|u|^2\,d\Vol_g\quad\text{for all }p\in B(0,9).$$
Therefore, by \eqref{eq:sse}, we have that
\begin{equation}\label{eq:ssev}
c_3\le\int_{B(p,1)}|v|^2\,dx\le c_4\quad\text{for all }p\in B(0,9).
\end{equation}
Here, the positive constants $c_3$ and $c_4$ depend only on $D_1$ and $D_2$ in \eqref{eq:sse} and $\M$.

Notice also that there are positive constants $c_5$ and $c_6$ depending only on $\M$ such that
\begin{equation}\label{eq:HuHv}
c_5r^{n-1}\,\HH^{n-1}(\NN_v\cap B(0,1))\le\HH^{n-1}(\NN_u\cap B(x,r))\le c_6r^{n-1}\,\HH^{n-1}(\NN_v\cap B(0,1)).
\end{equation}
To prove \eqref{eq:nodalinB}, it now suffices to prove that
\begin{prop}\label{prop:nodalinV}
Let $v$ be an eigenfunction with eigenvalue $\mu^2\gg1$ on $(B(0,10),\tilde g)$. Assume that \eqref{eq:ssev} holds. Then
\begin{equation}\label{eq:nodalinV}
E_1\mu\le\HH^{n-1}(\NN_v\cap B(0,1))\le E_2\mu
\end{equation}
for some positive constants $E_1$ and $E_2$ depending only on $c_3$ and $c_4$ in \eqref{eq:ssev} and $\M$. 
\end{prop}

Proposition \ref{prop:nodalinV} is proved following the estimates of the nodal sets in Donnelly-Fefferman \cite{DF}. We need to establish the growth estimate of the eigenfunction $v$. First, one has the local growth estimate of $v$. (See also Lin \cite{L}.) That is, for $0<\delta\le1$, there is a positive constant $c_7$ such that
$$\max_{B(p,2\delta)}|v|\le e^{c_7\mu}\max_{B(p,\delta)}|v|\quad\text{for all }p\in B(0,1).$$
We then need to determine the dependence of $c_7$ in the above inequality in a global setting (i.e. uniform for $p\in\M$ as $\mu\to\infty$.) As proved in \cite[Theorem 4.2 (ii)]{DF}, $c_7$ depends on
\begin{enumerate}[(i).]
\item the upper bound of the absolute values of the sectional curvatures of $(B(0,5),\tilde g)$,
\item the upper bound of the diameter of $(B(0,5),\tilde g)$.
\item $\|v\|_{L^\infty(B(0,5))}$.
\end{enumerate}
Here, the upper bound of the absolute values of the sectional curvatures of $B(0,5)$ is $K/2$, while the diameter $B(0,5)$ is uniformly bounded.

In (iii), $c_7$ depends on the global normalization that $\|v\|_{L^\infty(B(0,5))}=1$. In fact, this normalization is the only condition that is needed to get from the above local growth estimate to the global estimates of eigenfunctions e.g. vanishing order estimate and nodal set estimate. See the introduction in \cite[Section 1]{DF}.

In Proposition \ref{prop:nodalinV}, instead of the $L^\infty$ normalization, we see from \eqref{eq:ssev} that the 
\begin{equation}\label{eq:vL2}
\|v\|_{L^2(B(p,5))}\approx1\quad\text{for all }p\in B(0,1)
\end{equation} 
by a simple covering argument and the fact that $\|v\|_{L^2(B(x,1))}\approx1$ for all $x\in B(0,9)$ in \eqref{eq:ssev}. These two normalization obviously are not equivalent, that is, from H\"omander \cite{Ho}, generally one only has that
$$\|v\|_{L^2(B(0,5))}\lesssim\|v\|_{L^\infty(B(0,5))}\lesssim\mu^{(n-1)/2}\|v\|_{L^2(B(0,6))}.$$
However, in the proof of Proposition 4.1 and Theorem 4.2 in \cite[Section 4]{DF}, the only fact that one uses $\|v\|_{L^\infty(B(0,5))}=1$ is that 
\begin{equation}\label{eq:vp0}
v(p_0)\approx1\quad\text{for some }p_0\in B(0,5).
\end{equation} 
But it follows from \eqref{eq:vL2} trivially as well. That is, one has that
$$c_8\le v(p_0)\le c_9\quad\text{for some }p_0\in B(0,5),$$
in which $c_8$ and $c_9$ are positive constants that depend only on $c_3$, $c_4$ in \eqref{eq:ssev} and $K$. We thus arrive at the crucial global growth estimate of the eigenfunction $v$.
\begin{lemma}[Global growth estimate of eigenfunctions]\label{lemma:growth}
Let $v$ be defined in Proposition \ref{prop:nodalinV}.  Assume that \eqref{eq:ssev} holds. Then there is $0<\delta\le1$ such that
$$\max_{B(p,2\delta)}|v|\le e^{c_7\mu}\max_{B(p,\delta)}|v|\quad\text{for all }p\in B(0,1),$$
where $c_7$ depends only on $c_3$ and $c_4$ in \eqref{eq:ssev} and $\M$.
\end{lemma}

Once the above global growth estimate of the eigenfunction $v$ in $B(0,1)$ is established, we can proceed to prove Proposition \ref{prop:nodalinV}. This is done by studying the holomorphic extension of $v$ to a complex neighborhood of $B(0,1)$. 

Since the metric $\tilde g$ is analytic in $B(0,10)$, the eigenfunction $v$ extends to a complex region $\{z\in\C^n:|z|<10\}$. Let $z$ denote the complex variable and $x$ the real variable. Then
\begin{equation}\label{eq:zvx}
\sup_{|z|<\rho_1}|v(z)|\le e^{c_8\mu}\sup_{|x|<1}|v(x)|
\end{equation}
for $\rho_1<1$. Here, $c_8$ depends on $c_7$ in Lemma \ref{lemma:growth} and $K$. See \cite[Lemma 7.1]{DF}. Together with Lemma \ref{lemma:growth}, we have that the fundamental estimate to derive both the lower and upper bounds in Proposition \ref{prop:nodalinV}:
\begin{equation}\label{eq:growthinC}
\sup_{|z|<\rho_1}|v(z)|\le e^{c_9\mu}\sup_{|x|<\rho_2}|v(x)|
\end{equation}
for $\rho_2<\rho_1<1$. Here, $c_9$ depends on $c_8$ in Lemma \ref{lemma:growth} and $K$, which in term depends on $c_3$ and $c_4$ in \eqref{eq:ssev} and $\M$. See \cite[Lemma 7.2]{DF}.

We next derive the lower and upper bounds of the nodal set of the eigenfunction $v$ in $B(0,1)$, using \eqref{eq:growthinC}. Because the arguments follows the one in \cite[Section 7]{DF}, we only make clear of the dependence of these bounds on the necessary constants appeared in Proposition \ref{prop:nodalinV}. 

\subsection{Lower bounds}
One can derive from elliptic theory that $B(p,a_1\mu^{-1})$ contains a nodal point of $v$ if $a_1$ is large enough. If a ball $B(p,20a_1\mu^{-1})\subset B(0,1)$ satisfies that
\begin{equation}\label{eq:doubling}
\int_{B(p,20a_1\mu^{-1})}|v|^2\le a_2\int_{B(p,10a_1\mu^{-1})}|v|^2
\end{equation}
for some positive constant $a_2$ independent of $\mu$, then
$$\HH^{n-1}(\NN_v\cap B(p,10a_1\mu^{-1}))\ge a_3\mu^{-(n-1)}$$
for some positive constant $a_3$ independent of $\mu$. See \cite[Section 7]{DF}, also Colding-Minicozzi \cite{CM}. To estimate the lower bound the measure of $\NN_v$ in $B(0,1)$, cover $B(0,1)$ by balls $B(p,20a_1\mu^{-1})$ by choosing $B(p,a_1\mu^{-1})$ as the maximal disjoint family in $B(0,1)$. Then \eqref{eq:growthinC} guarantees that at least half of $B(p,20a_1\mu^{-1})$ satisfy \eqref{eq:doubling}, (which in \cite{CM} are called the ``good'' balls satisfying the doubling condition \eqref{eq:doubling}.) One then just needs to sum over the measure of the nodal set in these balls $B(p,20a_1\mu^{-1})$ and get the lower bound of 
$$\HH^{n-1}(\NN_v\cap B(0,1))\ge E_1\mu,$$
where $E_1$ depends on $c_9$ in \eqref{eq:growthinC} and the manifold $\M$. Notice that we can control the overlapping of the balls $B(p,20a_1\mu^{-1})$ of the chosen covering in a uniform way depending only on $\M$. (See also Section \ref{sec:enodal} for a precise control of the overlap for a similar cover.) Hence, $E_1$ depends on $c_3$ and $c_4$ in \eqref{eq:ssev} and $\M$, finishing the lower bound in Proposition \ref{prop:nodalinV}.

\subsection{Upper bounds}
The upper bound in Proposition \ref{prop:nodalinV} follows an upper bound of the measure of the nodal set of a holomorphic function in $\{z\in\C^n:|z|<1\}$ by \cite[Proposition 6.7]{DF}. There the upper bounds follows the growth estimate in the complex region \ref{eq:growthinC} directly. Then
$$\HH^{n-1}(\NN_v\cap B(0,1))\le E_2\mu,$$
where $E_2$ depends on $c_9$ in \eqref{eq:growthinC} and the manifold $\M$. Hence, $E_2$ depends on $c_3$ and $c_4$ in \eqref{eq:ssev} and $\M$, finishing the upper bound in Proposition \ref{prop:nodalinV}.

\section{Distribution of the nodal sets}\label{sec:enodal}
Let $\{u_k\}_{k=1}^\infty$ be a sequence of eigenfunctions that is \eqref{eq:compnodal} is true at a small scale $r(\lambda)$ defined in the introduction. Then $r_k=r(\lambda_k)\to0$ as $k\to\infty$. By Theorem \ref{thm:nodalinB}, we have that 
$$E_1\Vol(B(x,r_k)\le\frac{1}{\lambda_k}\HH^{n-1}(\NN_{u_k}\cap B(x,r_k))\le E_2\Vol(B(x,r_k))\quad\text{as }k\to\infty$$
uniformly for all $x\in\M$, where the positive constants $E_1$ and $E_2$ depend only on $D_1$ and $D_2$ in \eqref{eq:sse} and $\M$.

In this section, we use covering arguments to prove the comparability of the hypersurface volume form of $\{\NN_{u_k}\}_{k=1}^\infty$ in Theorem \ref{thm:nodalonM}.

Select a maximal disjoint family of balls $\{B(x_i,r_k/2)\}_{i=1}^N$ in $\M$. Then it follows immediately from the maximality that the family $\{B(x_i,r_k)\}_{i=1}^N$ is a cover of $\M$. Moreover, the overlapping of $\{B(x_i,r_k)\}_{i=1}^N$ at all $x\in\M$ is uniformly bounded. That is, for any $x\in\M$, suppose that $x$ is an element of $B(x_{i_1},r_k), \cdots, B(x_{i_L},r_k)$. Then $B(x_{i_l},r_k/2)\subset B(x,2r_k)$ for all $l=1,...,L$. But $B(x_{i_1},r_k/2), \cdots, B(x_{i_L},r_k/2)$ are disjoint, so
$$\sum_{l=1}^L\Vol(B(x_{i_l},r_k/2))\le\Vol(B(x,2r_k),$$
which implies that $L\le L_\M$ for some positive constant $L_\M$ depending only on $\M$. 

For a nonnegative function $f\in C(\M)$, if $f\equiv0$, then Theorem \ref{thm:nodalonM} holds trivially; if $f\not\equiv0$, then 
\begin{equation}\label{eq:fge0}
\int_\M f\,d\Vol>0.
\end{equation}
We compute the lower and upper bounds of $\int_{\NN_{u_k}}f\,d\HH^{n-1}$ as follows.

\subsection{Lower bounds}
Observe that
$$\sum_{i=1}^N\int_{\NN_{u_k}\cap B(x_i,r_k)}f\,d\HH^{n-1}\le L_\M\int_{\NN_{u_k}}f\,d\HH^{n-1},$$
since the overlapping of the family $\{B(x_i,r_k)\}_{i=1}^N$ at all $x\in\M$ is uniformly bounded by $L_\M$. Hence,
\begin{eqnarray*}
&&\int_{\NN_{u_k}}f\,d\HH^{n-1}\\
&\ge&\frac{1}{L_\M}\sum_{i=1}^N\int_{\NN_{u_k}\cap B(x_i,r_k)}f\,d\HH^{n-1}\\
&=&\frac{1}{L_\M}\sum_{i=1}^Nf(\bar x_i)\HH^{n-1}(\NN_{u_k}\cap B(x_i,r_k))\quad\text{for some }\bar x_i\in B(x_i,r_k)\\
&\ge&\frac{E_1}{L_\M}\cdot\lambda_k\sum_{i=1}^Nf(\bar x_i)\Vol(B(x_i,r_k))\quad\text{by Theorem \ref{thm:nodalinB}}\\
&\ge&\frac{E_1}{L_\M}\cdot\lambda_k\sum_{i=1}^N\left(\int_{B(x_i,r_k)}f(x)\,d\Vol-\int_{B(x_i,r_k)}|f(x)-f(\bar x_i)|\,d\Vol\right)\\
&\ge&\frac{E_1}{L_\M}\cdot\lambda_k\left(\sum_{i=1}^N\int_{B(x_i,r_k)}f(x)\,d\Vol-\max_{d(x,y)\le r_k}|f(x)-f(y)|\sum_{i=1}^N\Vol(B(x_i,r_k))\right)\\
&\ge&\frac{E_1}{L_\M}\cdot\lambda_k\left(\int_\M f\,d\Vol-\max_{d(x,y)\le r_k}|f(x)-f(y)|\cdot\Vol(\M)\right)\\
&\ge&\frac{E_1}{2L_\M}\cdot\lambda_k\int_\M f\,d\Vol\quad\text{as }k\to\infty,
\end{eqnarray*}
since $\int_\M f\,d\Vol>0$ by \eqref{eq:fge0}. Here, we use the fact that $r_k\to0$ as $k\to\infty$ and $f\in C(\M)$ is uniformly continuous on $\M$.

\subsection{Upper bounds}
Compute that
\begin{eqnarray*}
&&\int_{\NN_{u_k}}f\,d\HH^{n-1}\\
&\le&\sum_{i=1}^N\int_{\NN_{u_k}\cap B(x_i,r_k)}f\,d\HH^{n-1}\\
&=&\sum_{i=1}^Nf(\bar x_i)\HH^{n-1}(\NN_{u_k}\cap B(x_i,r_k))\quad\text{for some }\bar x_i\in B(x_i,r_k)\\
&\le&E_2\cdot\lambda_k\sum_{i=1}^Nf(\bar x_i)\Vol(B(x_i,r_k))\quad\text{by Theorem \ref{thm:nodalinB}}\\
&\le&E_2\cdot\lambda_k\sum_{i=1}^N\left(\int_{B(x_i,r_k)}f(x)\,d\Vol+\int_{B(x_i,r_k)}|f(x)-f(\bar x_i)|\,d\Vol\right)\\
&\le&E_2\cdot\lambda_k\left(\sum_{i=1}^N\int_{B(x_i,r_k)}f(x)\,d\Vol+\max_{d(x,y)\le r_k}|f(x)-f(y)|\sum_{i=1}^N\Vol(B(x_i,r_k))\right)\\
&\le&E_2\cdot\lambda_k\left(L_\M\int_\M f\,d\Vol+\max_{d(x,y)\le r_k}|f(x)-f(y)|\cdot L_\M\Vol(\M)\right)\\
&\le&2L_\M E_2\cdot\lambda_k\int_\M f\,d\Vol\quad\text{as }k\to\infty,
\end{eqnarray*}
since $\int_\M f\,d\Vol>0$ by \eqref{eq:fge0}. Here, we use the fact that $r_k\to0$ as $k\to\infty$ and $f\in C(\M)$ is uniformly continuous on $\M$.

\section*{Acknowledgments}
I would like to thank Steve Zelditch for reading the note and offering suggestions that helped to improve the presentation.

\end{document}